\documentclass[10pt,twoside]{article}
\usepackage{graphicx}
\usepackage{amsmath}
\usepackage{Latex-document}

\title{\bf Index Iteration Theory for Symplectic \vskip -2mm Paths with
Applications to Nonlinear \vskip -2mm Hamiltonian Systems\vskip
6mm}
\author{Yiming Long\vspace*{-0.5cm}\thanks{Nankai Institute of
Mathematics, Nankai University, Tianjin 300071, China. E-mail:
longym@nankai.edu.cn}}
\date{\vspace{-8mm}}

\begin{document}

\maketitle

\thispagestyle{first} \setcounter{page}{303}

\begin{abstract}

\vskip 3mm

In recent years, we have established the iteration theory of
the index for symplectic matrix paths and applied it to periodic
solution problems of nonlinear Hamiltonian systems. This paper
is a survey on these results.

\vskip 4.5mm


\end{abstract}

\vskip 12mm

Since P. Rabinowitz's pioneering work \cite{Rab1} of 1978,
variational methods have been widely used in the study of
existence of solutions of Hamiltonian systems. But how to study
the geometric multiplicity and stability of periodic solution
orbits obtained by variational methods has kept to be a difficulty
problem. For example let $x=x(t)$ be a $\tau$-periodic solution of
a Hamiltonian system
\begin{equation}
    \dot{x}(t) = JH'(x(t)),  \quad \forall t\in {\bf R}. \label{1}
\end{equation}
The $m$-th iteration $x^m$ of $x$ is defined by induction $m-1$ times
via $x(t+\tau)=x(t)$ for $t>0$. It runs $m$-times along the orbit
of $x$. Geometrically these iterations produce the same solution orbit
of (\ref{1}), but they are different as critical points of corresponding
functionals. This multiple covering phenomenon causes major difficulties
in the study.

A natural way to study solution orbits found by variational
methods is to study the Morse-type index sequences of their
iterations. But when one studies general Hamiltonian systems, the
Morse indices of the critical points of the corresponding
functional are always infinite. To overcome this difficulty, in
their celebrated paper \cite{CoZ1} of 1984, C. Conley and E.
Zehnder defined an index theory for any non-degenerate paths in
${\rm Sp}(2n)$ with $n\geq 2$, i.e., the so called Conley-Zehnder
index theory. This index theory was further defined for
non-degenerate paths in ${\rm Sp}(2)$ by E. Zehnder and the author
in \cite{LZe1} of 1990. The index theory for degenerate linear
Hamiltonian systems was defined by C. Viterbo in \cite{Vit1} and
the author in \cite{Lon1} of 1990 independently. In \cite{Lon6} of
1997, this index was extended to any degenerate symplectic matrix
paths.

Motivated by the iteration theories for the Morse type index theories
established by R. Bott in 1956 and by I. Ekeland in 1980s, in
recent years the author extended the index theory mentioned above,
introduced an index function theory for symplectic matrix paths, and
established the iteration theory for the index theory of
symplectic paths. Applying this index iteration theory to nonlinear
Hamiltonian systems, interesting results on periodic solution problems of
Hamiltonian systems are obtained. Here a brief survey is given on these
subjects. Readers are referred to the author's recent book \cite{Lon11}
for further details.

\section{Index function theory for symplectic paths} \label{section 1}
\setzero

\vskip-5mm \hspace{5mm}

As usual we define the symplectic group by ${\rm Sp}(2n) = \{M\in {\rm GL}({\bf R}^{2n})\,|\,M^TJM$ $=J\}$, where
$J=\left(\begin{matrix} 0&-I\cr
                       I&0\cr \end{matrix}\right)$, $I$ is the identity matrix
on ${\bf R}^n$, and $M^T$ denotes the transpose of $M$. For
$\omega\in {\bf U}$, the unit circle in the complex plane ${\bf C}$,
we define the $\omega$-singular subset in ${\rm Sp}(2n)$ by
${\rm Sp}(2n)_{\omega}^0 = \{M\in {\rm Sp}(2n)\,|\,
\omega^{-n}\det(\gamma(\tau)-\omega I)=0\}$.
Here for any $M\in {\rm Sp}(2n)_{\omega}^0$, we define the orientation
of ${\rm Sp}(2n)_{\omega}^0$ at $M$ by the positive direction
$\frac{d}{dt}M\exp(tJ)|_{t=0}$. Since the fundamental solution of a
general linear Hamiltonian system with continuous symmetric periodic
coefficient $2n\times 2n$ matrix function $B(t)$,
\begin{equation}
    \dot{x}(t) = JB(t)x(t),  \quad \forall t\in {\bf R}, \label{2}
\end{equation}
is a path in ${\rm Sp}(2n)$ starting from the identity, for $\tau>0$
we define the set of symplectic matrix paths by
${\cal P}_{\tau}(2n) = \{\gamma\in C([0,\tau],{\rm Sp}(2n))\,|\,
\gamma(0)=I\}$. For any two path $\xi$ and $\eta:[0,\tau]\to {\rm Sp}(2n)$
with $\xi(\tau)=\eta(0)$, as usual we define $\eta\ast\xi(t)$ by
$\xi(2t)$ if $0\le t\le\tau/2$, and $\eta(2t-\tau)$ if
$\tau/2\le t\le\tau$. We define a special path
$\zeta:[0,\tau]\to{\rm Sp}(2n)$ by
$$ \zeta(t) =
{\rm diag}(2-\frac{t}{\tau},\ldots, 2-\frac{t}{\tau},
               (2-\frac{t}{\tau})^{-1},\ldots, (2-\frac{t}{\tau})^{-1}),
                 \quad {\rm for}\;0\le t\le \tau.  $$

  {\bf Definition 1.} (cf. \cite{Lon8}) {\it For any $\tau>0$,
$\omega\in {\bf U}$, and $\gamma\in {\cal P}_{\tau}(2n)$,
we define
the $\omega$-nullity of $\gamma$ by
\begin{equation}
\nu_{\omega}(\gamma) =
\dim_{\bf C}\ker_{\bf C}(\gamma(\tau)-\omega I). \label{3}
\end{equation}
If $\gamma$ is $\omega$ non-degenerate, i.e., $\nu_{\omega}(\gamma)=0$, we
define the $\omega$-index of $\gamma$ by the intersection number
\begin{equation}
i_{\omega}(\gamma) = [{\rm Sp}(2n)_{\omega}^0:\gamma\ast\zeta].
\label{4}
\end{equation}
If $\gamma$ is $\omega$ degenerate, i.e., $\nu_{\omega}(\gamma)>0$,
we let ${\cal F}(\gamma)$ be the set of all open neighborhoods of $\gamma$
in ${\cal P}_{\tau}(2n)$, and define
\begin{equation}
i_{\omega}(\gamma) = \sup_{U\in{\cal F}(\gamma)}\inf\{i_{\omega}(\beta)\,|\,
          \beta\in U, \,\nu_{\omega}(\beta)=0\}. \label{5}
\end{equation}
Then we call $(i_{\omega}(\gamma), \nu_{\omega}(\gamma))\in {\bf
Z}\times \{0,1,\ldots,2n\}$ the index function of $\gamma$ at
$\omega$.}

  The relation of this index $(i_1(\gamma), \nu_1(\gamma))$ with the
Morse index of $\tau$-periodic solutions of the problem (1) was
proved by C. Conley, E. Zehnder, and the author in \cite{CoZ1},
\cite{LZe1}, and \cite{Lon1} (cf. Theorem 6.1.1 of \cite{Lon11}).

\section{\hskip -2mm Iteration theory of the index for symplectic paths}
\label{section 2}\setzero\vskip-5mm \hspace{5mm }

  Given a path $\gamma\in {\cal P}_{\tau}(2n)$, its iteration is defined
inductively by $\gamma(t+\tau)=\gamma(t)\gamma(\tau)$ for $t\ge 0$, i.e.,
\begin{equation}
\gamma^m(t) = \gamma(t-j\tau)\gamma(\tau)^j, \qquad j\tau\le t\le (j+1)\tau,
        j=0,1,\ldots, m-1,  \label{6}
\end{equation}
for any $m$ in the natural integer set ${\bf N}$. For our applications of
this index theory to nonlinear Hamiltonian systems, we are facing two
types of problems:

  $\langle 1\rangle$ knowing the end point $\gamma(\tau)$ of a path
$\gamma\in {\cal P}_{\tau}(2n)$, the initial index
$(i_1(\gamma),\nu_1(\gamma))$,
and the iteration time $m$, want to find the index $i_1(\gamma^m)$ of the
$m$-th iterated path $\gamma^m$;

  $\langle 2\rangle$ knowing the end point $\gamma(\tau)$ of a path
$\gamma\in {\cal P}_{\tau}(2n)$, the initial index
$(i_1(\gamma),\nu_1(\gamma))$, and the index
$(i_1(\gamma^m),\nu_1(\gamma^m))$
of the $m$-th iterated path $\gamma^m$, want to find the iteration time $m$.

  To solve these problems, we first generalize Bott's formula of the iterated
Morse index for closed geodesics to the index theory for general symplectic
paths:

{\bf Theorem 2} (cf. \cite{Lon8}){\bf .} {\it For any $\tau>0$,
$\gamma\in {\cal P}_{\tau}(2n)$, $z\in{\bf U}$, and $m\in {\bf
N}$, there hold: }
\begin{equation}
  i_z(\gamma^m) = \sum_{\omega^m=z}i_{\omega}(\gamma), \quad
  \nu_z(\gamma^m) = \sum_{\omega^m=z}\nu_{\omega}(\gamma). \label{7}
\end{equation}

By (\ref{7}) it is easy to see that the mean index
$\hat{i}(\gamma) = \lim_{m\to +\infty}i_1(\gamma^m)/m$ for any
$\gamma\in {\cal P}_{\tau}(2n)$ is always a finite real number.

  To further solve the problems $\langle 1\rangle$ and
$\langle 2\rangle$, we need to go beyond the Bott-type formula
(\ref{7}). For a given path $\gamma$ we consider to deform
it to a new path $\eta$ in ${\cal P}_{\tau}(2n)$ so that
\begin{equation}
i_1(\gamma^m)=i_1(\eta^m),\quad \nu_1(\gamma^m)=\nu_1(\eta^m), \quad
         \forall m\in {\bf N}, \label{8}
\end{equation}
and that $(i_1(\eta^m),\nu_1(\eta^m))$ is easy enough to compute. This
leads to finding homotopies $\delta:[0,1]\times[0,\tau]\to {\rm Sp}(2n)$
starting from $\gamma$ in ${\cal P}_{\tau}(2n)$ and keeping the end
points of the homotopy always stay in a certain suitably chosen maximal
subset of ${\rm Sp}(2n)$ so that (\ref{8}) always holds. By (\ref{7}), this
set is defined to be the path connected component
$\Omega^0(M)$ containing $M=\gamma(\tau)$ of the set
\begin{eqnarray}
  \Omega(M)=\{N\in{\rm Sp}(2n)\,&|&\,\sigma(N)\cap{\bf U}=\sigma(M)\cap{\bf U},\,
{\rm and}\;  \nonumber\\
 &&\qquad \nu_{\lambda}(N)=\nu_{\lambda}(M)\;\forall\,
\lambda\in\sigma(M)\cap{\bf U}\}. \label{9}
\end{eqnarray}
Here we call $\Omega^0(M)$ the {\it homotopy component} of $M$ in
${\rm Sp}(2n)$.

Using normal forms of symplectic matrices (cf. \cite{LoD1}, \cite{HaL1}),
we then decompose $\gamma(\tau)$ within $\Omega^0(\gamma(\tau))$ into
product of 10 special $2\times 2$ and $4\times 4$ symplectic normal form
matrices, which we call {\it basic normal forms}. Correspondingly by the
homotopy invariance and symplectic additivity of the index theory, the
computations in (\ref{8}) are reduced to iterations of those paths in
${\rm Sp}(2)$ or ${\rm Sp}(4)$ whose end points are one of the 10 basic
normal form matrices. The study of the index for iterations of any
symplectic paths is carried out for paths in ${\rm Sp}(2)$ via the
${\bf R}^3$-cylindrical coordinate representation of ${\rm Sp}(2)$, then
for hyperbolic and elliptic paths in ${\rm Sp}(2n)$. This yields the precise
iteration formula obtained in \cite{Lon10} of the index theory for any
symplectic path $\gamma\in {\cal P}_{\tau}(2n)$ in terms of the basic
norm form decomposition of $\gamma(\tau)$, $(i(\gamma,1),\nu(\gamma,1))$,
and the iteration time $m$.

  For any $M\in {\rm Sp}(2n)$, its splitting numbers at an
$\omega\in {\bf U}$ is defined in \cite{Lon8} by
\begin{equation}
S_M^{\pm}(\omega) =
\lim_{\epsilon\to 0^+}i_{\omega\exp(\pm\sqrt{-1}\epsilon)}(\gamma)
- i_{\omega}(\gamma),  \label{10}
\end{equation}
via any $\gamma\in{\cal P}_{\tau}(2n)$ satisfying $\gamma(\tau)=M$. Then
it is proved that the splitting numbers of $M$ at $\omega$
can be characterized algebraically.

  Motivated by the precise iteration formulae of \cite{Lon10}, the
following second index iteration formula of any symplectic path is
established by C. Zhu and the author. Here we denote by
$(i(\gamma,m),\nu(\gamma,m))$ $=$ $(i_1(\gamma^m),\nu_1(\gamma^m))$.

  {\bf Theorem 3} (cf. \cite{LZh1}){\bf .} {\it For any $\tau>0$,
$\gamma\in {\cal P}_{\tau}(2n)$, and $m\in {\bf N}$,
there holds:
\begin{eqnarray}
i(\gamma,m) &=& m(i(\gamma,1) + S^+_M(1)- C(M))   \nonumber\\
&& + 2\sum_{\theta\in(0,2\pi)}E(\frac{m\theta}{2\pi})
      S^-_M(e^{\sqrt{-1}\theta}) - (S^+_M(1)+C(M)), \label{11}
\end{eqnarray}
where $M=\gamma(\tau)$,
$C(M)=\sum_{0<\theta<2\pi}S_M^-(e^{\sqrt{-1}\theta})$, and
$E(a)=\min\{k\in {\bf Z}\,|\,k\ge a\}$ for any $a\in {\bf R}$. }

In order to solve problems on nonlinear Hamiltonian systems, various
index iteration inequalities for any path $\gamma\in{\cal P}_{\tau}(2n)$
and $m\in {\bf N}$ are proved by D. Dong, C. Liu, C. Zhu and the author
in \cite{DoL1}, \cite{LLo1}, \cite{LLo2}, and \cite{LZh1}.

{\bf Theorem 4.} {\it For any $\gamma\in{\cal P}_{\tau}(2n)$ and
$m\in{\bf N}$, the following iteration inequalities always hold.

  {\sl Estimate via mean index} (cf. \cite{LLo1}, \cite{LLo2}):
\begin{equation}
m\hat{i}(\gamma)-n \le i(\gamma,m) \le
m\hat{i}(\gamma)+n-\nu(\gamma,m).   \label{12}
\end{equation}

  {\sl Estimate via initial index} (cf. \cite{LLo3}):
\begin{eqnarray}
m(i(\gamma,1) &+& \nu(\gamma,1) - n) + n - \nu(\gamma,1) \;\le\;
i(\gamma,m)      \nonumber\\
&\le & m(i(\gamma,1)+n)-n - (\nu(\gamma,m) - \nu(\gamma,1)). \label{13}
\end{eqnarray}

  {\sl Successive index estimate} (cf. \cite{LZh1}):
\begin{eqnarray}
 \nu(\gamma,m) - \frac{e(\gamma(\tau))}{2}
&\le& i(\gamma,m+1) - i(\gamma,m) - i(\gamma,1)  \nonumber\\
&\le& \nu(\gamma,1) -
\nu(\gamma,m+1) + \frac{e(\gamma(\tau))}{2}. \label{14}
\end{eqnarray}
Here we define $e(M)$ to be the total multiplicity of eigenvalues of
$M$ on ${\bf U}$ and call it the {\it elliptic height} of $M$. }

  A consequence of the iteration inequality (\ref{13}) together with
the necessary and sufficient conditions for any equality in (\ref{13})
to hold for some $m$ yields a new proof of the following theorem of
D. Dong and the author on controlling the iteration time
$m$ via indices:

{\bf Theorem 5} (cf. \cite{DoL1}){\bf .} {\it For any $\gamma\in
{\cal P}_{\tau}(2n)$ and $m\in{\bf N}$, suppose $i(\gamma,m)\le
n+1$, $i(\gamma,1)\ge n$, and $\nu(\gamma,1)\ge 1$. Then $m=1$. }

Note also that the inequality (\ref{14}) yields a way to estimate the
ellipticity of solution orbits of Hamiltonian systems obtained by
variational methods via their iterated indices.

  In order to study the properties of solution orbits of the system
(\ref{1}) on a given energy hypersurface, when the number of orbits
is finite, we need to study common properties of any given finite
family of symplectic paths $\gamma_j\in {\cal P}_{\tau_j}(2n)$ with
$1\le j\le q$. This leads to the following common index jump theorem
of C. Zhu and the author proved in \cite{LZh1}. For any
$\gamma\in{\cal P}_{\tau}(2n)$, its $m$-th index jump
${\cal G}_m(\gamma)$ is defined to be the open interval
${\cal G}_m(\gamma) =(i(\gamma,m) + \nu(\gamma,m) - 1, i(\gamma,m+2))$.

{\bf Theorem 6} (cf. \cite{LZh1}){\bf .} {\it Let $\gamma_j\in
{\cal P}_{\tau_j}(2n)$ with $1\le j\le q$ satisfying
\begin{equation}
\hat{i}(\gamma_j)>0, \quad i(\gamma_j,1)\ge n, \qquad 1\le j\le q.
\label{15}
\end{equation}
Then there exist infinitely many positive integer tuples
$(N,m_1,\ldots,m_q)\in {\bf N}^{q+1}$ such that
\begin{equation}
\emptyset\not= [2N-\kappa_1,2N+\kappa_2] \subset
      \bigcap_{j=1}^q{\cal G}_{2m_j-1}(\gamma_j),     \label{16}
\end{equation}
where
$\kappa_1=\min_{1\le j\le q}(i(\gamma_j,1)
+ 2S_{\gamma_j(\tau_j)}^+(1)- \nu(\gamma_j,1))$ and
$\kappa_2=\min_{1\le j\le q}i(\gamma_j,1)-1$. }

  In order to prove this theorem, we need to make each index jump to
be as big as possible, and to make their largest sizes happen
simultaneously to guarantee the existence of a non-empty largest
common intersection interval among them. By the term
$E(\frac{m\theta}{2\pi})$ in the abstract iteration formula
(\ref{11}), such a problem is reduced to a dynamical system problem
on a torus, and is solved by properties of closed additive subgroups
of tori.

\section{Applications to nonlinear Hamiltonian systems}\label{section 3}
\setzero\vskip-5mm \hspace{5mm }

So far, we have applied our index iteration theory to three important
problems on periodic solutions of nonlinear Hamiltonian systems. Let
$T>0$ and suppose $x$ is a non-constant $T$-periodic solution of the
nonlinear Hamiltonian system (\ref{1}). Suppose the minimal period of
$x$ is $\tau=T/k$ for some $k\in {\bf N}$. We denote by
$\gamma_x\in{\cal P}_{\tau}(2n)$ the fundamental solution of the
linearized Hamiltonian system (\ref{2}) at $x$ with $B(t)=H''(x(t))$,
and the iterated index of $x$ by
$(i(x,m),\nu(x,m)) = (i(\gamma_x,m),\nu(\gamma_x,m))$ for all
$m\in {\bf N}$.

\subsection{Prescribed minimal period solution problem}\label{subsection 3.1}
\setzero\vskip-4mm \hspace{4mm }

  In \cite{Rab1} of 1978, P. Rabinowitz posed a conjecture on whether the
Hamiltonian system possesses periodic solutions with prescribed minimal
period when the Hamiltonian function satisfies his superquadratic
conditions. This conjecture is studied by D. Dong and the author as an
application of our index iteration theory. Note that for a non-constant
$\tau$-periodic solution $x$ of the autonomous system (\ref{1}), the
condition on the nullity in Theorem 5 always holds. Thus Theorem 5 yields:

{\bf Theorem 7} (cf. \cite{DoL1}){\bf .} {\it For any non-constant
$\tau$-periodic solution $x$ of (\ref{1}), denote its minimal
period by $\tau/m$ for some $m\in{\bf N}$. Suppose
$i(x|_{[0,\tau]},1)\le n+1$ and $n\le i(x|_{[0,\tau/m]},1)$. Then
$m=1$, i.e., $\tau$ is the minimal period of $x$. }

  Here the first estimate on the index holds if $x$ is obtained by
minimax or minimization methods, and the second estimate on the index
holds if the Hamiltonian function $H$ is convex in a certain weak sense
along the orbit of $x$. This result reveals the intrinsic relationship
between the minimal period of a periodic solution and its indices, and unifies
all the results on Rabinowitz's conjecture under various convexity
conditions. Specially, it recovers the famous theorem of I. Ekeland
and H. Hofer in 1985 (cf. \cite{EkH1}) who solved Rabinowitz's conjecture for
convex superquadratic Hamiltonian systems.

\subsection{Periodic points of the Poincar\'e map of Lagrangian systems on
tori}\label{subsection 3.2}\setzero\vskip-4mm \hspace{4mm }

  In 1984, C. Conley stated a conjecture on whether the Poincar\'e map
of any $1$-periodic time dependent Hamiltonian system defined on the
standard torus $T^{2n}={\bf R}^{2n}/{\bf Z}^{2n}$ always possesses infinitely
many periodic points which are produced by contractible periodic solutions of
the corresponding Hamiltonian system on $T^{2n}$. A celebrated partial answer to
this conjecture was given by D. Salamon and E. Zehnder in 1992 (cf. \cite{SaZ1})
for a large class of symplectic manifolds on which every contractible integer
periodic solution of the Hamiltonian system has at least one Floquet multiplier
not equal to $1$. So far Conley conjecture is still open and seems far from
being completely understood.

  In \cite{Lon9}, we studied the Lagrangian system version of this conjecture.
Consider
\begin{equation}
{d\over{dt}}L_{\dot{x}}(t,x,\dot{x})-L_x(t,x,\dot{x}) = 0,
      \quad x\in {\bf R}^n,  \label{17}
\end{equation}
where $L_{\dot{x}}$ and $L_x$ denote the gradients of $L$ with respect
to $\dot{x}$ and $x$ respectively. The main result is the following:

{\bf Theorem 8} (cf. \cite{Lon9}){\bf .} {\it Suppose the
Lagrangian function $L$ satisfies

(L1) $L(t,x,p)={1\over 2}A(t)p\cdot p+V(t,x)$, where
${1\over 2}A(t)p\cdot p\geq \lambda |p|^2$ for all
$(t,p)\in {\bf R}\times {\bf R}^n$ and some fixed constant $\lambda>0$.

(L2) $A\in C^3({\bf R},{\cal L}_s({\bf R}^n))$, $V\in C^3({\bf
R}\times{\bf R}^n,{\bf R})$, both $A$ and $V$ are $1$-periodic in all
of their variables, where ${\cal L}_s({\bf R}^n)$ denotes the set of
$n\times n$ real symmetric matrices.

Then the Poincar\'e map $\Psi$ of the system (\ref{17}) possesses
infinitely many periodic points on $TT^n$ produced by contractible
integer periodic solutions of the system (\ref{17}). }

  In the proof of Theorem 8, the above inequality (\ref{12}) plays a
crucial role. By this inequality, at very high iteration level, a global
homological injection map can be constructed which maps a generator of
a certain non-trivial local critical group to a nontrivial homology
class $[\sigma]$ in a global homology group, if the number of
contractible integer periodic solution towers of the system (\ref{17})
is finite. But on the other hand, by a technique of V. Bangert
and W. Klingenberg in \cite{BaK1}, it is shown that this homology
class $[\sigma]$ must be trivial globally. This contradiction then
yields the conclusion of Theorem 8.

\subsection{Closed characteristics on convex compact hypersurfaces}
\label{subsection 3.3}\setzero\vskip-4mm \hspace{4mm }

  Denote the set of all compact strictly convex $C^2$-hypersurfaces
in ${\bf R}^{2n}$ by ${\cal H}(2n)$. For $\Sigma\in{\cal H}(2n)$ and
$x\in\Sigma$, let $N_{\Sigma}(x)$ be the outward normal unit vector
at $x$ of $\Sigma$. We consider the problem of finding $\tau>0$ and a
curve $x\in C^1([0,\tau],{\bf R}^{2n})$ such that
\begin{equation}
\left\{\begin{array}{ccl}
  \dot{x}(t) &=& JN_{\Sigma}(x(t)), \quad x(t)\in\Sigma,
                       \quad\forall t\in {\bf R}, \\
  x(\tau) &=& x(0). \end{array}\right. \label{18}
\end{equation}
A solution $(\tau,x)$ of the problem (\ref{18}) is called a {\it closed
characteristic} on $\Sigma$. Two closed characteristics $(\tau,x)$ and
$(\sigma,y)$ are {\it geometrically distinct}, if
$x({\bf R})\not= y({\bf R})$. We denote by ${\cal T}(\Sigma)$ the set of
all geometrically distinct closed characteristics $(\tau,x)$ on $\Sigma$
with $\tau$ being the minimal period of $x$. Note that the problem
(\ref{18}) can be described in a Hamiltonian system version and solved
by variational methods. A closed characteristic $(\tau,x)$ is
{\it non-degenerate}, if $1$ is a Floquet multiplier of $x$ of precisely
algebraic multiplicity $2$, and is {\it elliptic}, if all the Floquet
multipliers of $x$ are on ${\bf U}$. Let $^{\#}A$ denote the total number
of elements in a set $A$.

   This problem has been studied for more than 100 years since at
least A. M. Liapunov in 1892. A long standing conjecture on the
multiplicity of closed characteristics is whether
\begin{equation}
 ^{\#}\tilde{{\cal J}}(\Sigma)\ge n, \qquad \forall \,
\Sigma\in {\cal H}(2n). \label{19}
\end{equation}
The first break through on this problem in the global sense was
made by P. Rabinowitz \cite{Rab1} and A. Weinstein \cite{Wei1} in
1978. They proved $\,^{\#}{\cal T}(\Sigma)\ge 1$ for all
$\Sigma\in {\cal H}(2n)$. Besides many results under pinching
conditions, in 1987--1988, I. Ekeland-L. Lassoued, I. Ekeland-H.
Hofer, and A, Szulkin proved $\,^{\#}{\cal T}(\Sigma)\ge 2$ for
all $\Sigma\in {\cal H}(2n)$ and $n\ge 2$. In 1998, H. Hofer, K.
Wysocki, and E. Zehnder proved in \cite{HWZ1}: $\,^{\#}{\cal
T}(\Sigma)= 2$ or $+\infty$ for every $\Sigma\in{\cal H}(4)$. In
recent years C. Liu, C. Zhu, and the author gave the following
answers to the conjecture (\ref{19}):

{\bf Theorem 9} (cf. \cite{LZh1}){\bf .} {\it There holds
\begin{equation}
 ^{\#}{\cal T}(\Sigma)\ge [\frac{n}{2}]+1, \quad
\forall \Sigma\in {\cal H}(2n),  \label{20}
\end{equation}
where $[a]=\max\{k\in{\bf Z}\,|\,k\le a\}$ for any $a\in{\bf R}$.
Moreover, if all the closed characteristics on $\Sigma$ are non-degenerate,
then $\,^{\#}{\cal T}(\Sigma)\ge n$. }

{\bf Theorem 10} (cf. \cite{LLZ1}){\bf .} {\it For any $\Sigma\in
{\cal H}(2n)$, if $\Sigma$ is symmetric with respect to the
origin, i.e., $x\in\Sigma$ implies $-x\in\Sigma$, then
$\,^{\#}{\cal T}(\Sigma)\ge n$. }

  Very recently, Y. Dong and the author further proved the following result.

{\bf Theorem 11} (cf. \cite{DyL1}){\bf .} {\it Let $\Sigma\in
{\cal H}(2n)$ be $P$-symmetric with respect to the origin, i.e.,
$x\in\Sigma$ implies $Px\in\Sigma$, where $P={\rm
diag}(-I_{n-k},I_k,-I_{n-k},I_k)$ for some fixed integer $k\in
[0,n-1]$. Let $\Sigma(k)= \{(x,y)\in ({\bf
R}^k)^2\,|\,(0,x,0,y)\in\Sigma\}$. Suppose $\,^{\#}{\cal
T}(\Sigma(k)) \le k$ or $\,^{\#}{\cal T}(\Sigma(k))=+\infty$
holds. Then $\,^{\#}{\cal T}(\Sigma) \ge n-2k$. }

Proof of Theorem 11 depends on a new index iteration theory for
symplectic paths iterated by the formula
$\gamma(t+\tau)=P\gamma(t)P\gamma(\tau)$ for $t\ge 0$.

The second long standing conjecture on closed characteristics is
whether there always exists at least an elliptic closed
characteristic on any $\Sigma\in {\cal H}(2n)$. Up to the author's
knowledge, the existence of one elliptic closed characteristic on
$\Sigma\in {\cal H}(2n)$ was proved by I. Ekeland in 1990 when
$\Sigma$ is $\sqrt{2}$-pinched by two spheres, and by G.-F.
Dell'Antonio, B. D'Onofrio, and I. Ekeland in 1992 when
$\Sigma$ is symmetric with respect to the origin. Recently using an
enhanced version of the iteration estimate (\ref{14}) on the elliptic
height, based on results in \cite{Lon10} the following result was
further proved by C. Zhu and the author.

{\bf Theorem 12} (cf. \cite{LZh1}){\bf .} {\it For $\Sigma\in
{\cal H}(2n)$, suppose $\,^{\#}{\cal T}(\Sigma)<+\infty$. Then
there exists at least an elliptic closed characteristic on
$\Sigma$. Moreover, suppose $n\ge 2$ and $\,^{\#}{\cal
T}(\Sigma)\le 2[n/2]$. Then there exist at least two elliptic
elements in ${\cal T}(\Sigma)$. }

  The main ingredient in the proofs of Theorems 9 to 12 is our index
iteration theory mentioned above. To illustrate this method, we briefly
describe below the main idea in the proof of (\ref{20}) in Theorem 9.
Because each closed characteristic on $\Sigma$ corresponds to
infinitely many critical values of the related dual action functional,
our way to solve the problem is to study how the index intervals of
iterated closed characteristics cover the set of integers $2{\bf N}-2+n$
to count the number of closed characteristics on $\Sigma$. Suppose
$q = \,^{\#}\tilde{\cal J}(\Sigma)<+\infty$. In the proof
of the multiplicity claim (\ref{20}) of Theorem 9, the most important
ingredient is the following estimates:
\begin{eqnarray}
q
&\ge& ^{\#}\left((2{\bf N}-2+n)\cap
   \cap_{j=1}^q{\cal G}_{2m_j-1}(\gamma_{x_j})\right)  \nonumber\\
&\ge& ^{\#}\left((2{\bf N}-2+n)\cap
   [2N-\kappa_1,2N+\kappa_2]\right) \nonumber\\
&\ge& [\frac{n}{2}] + 1, \label{21}
\end{eqnarray}
The first inequality in (\ref{21}) is a new version of the
Liusternik-Schnirelman theoretical argument at the iterated index
level, which distinguishes solution orbits geometrically instead of
critical points only as usual methods do. The second inequality in
(\ref{21}) uses the common index jump Theorem 6. The last inequality
in (\ref{21}) uses the Morse theoretical approach. Roughly speaking,
the common index jump theorem picks up as many as possible points of
$2{\bf N}-2+n$ in the interval $[2N-\kappa_1,2N+\kappa_2]$ $\subset$
$\cap_{j=1}^q{\cal G}_{2m_j-1}(\gamma_{x_j})$, which yields a lower
bound for $^{\#}{\cal T}(\Sigma)$.

  As usual, a hypersurface $\Sigma\subset {\bf R}^{2n}$ is
star-shaped if the tangent hyperplane at any $x\in\Sigma$ does not
intersect the origin. Closed characteristics on $\Sigma$ can be
defined by (\ref{18}) too. In this case, the result
$\,^{\#}\tilde{\cal J}(\Sigma) \ge 1$ was proved by P. Rabinowitz in
\cite{Rab1} of 1978. Then multiplicity results were proved under
certain pinching conditions on star-shaped $\Sigma$. Recently, the
following result for the free case was proved by X. Hu and the
author:

{\bf Theorem 13} (cf. \cite{HuL1}){\bf .} {\it Let $\Sigma$ be a
star-shaped compact $C^2$-hypersurface in ${\bf R}^{2n}$. Suppose
all the closed characteristics on $\Sigma$ and all of their
iterates are non-degenerate. Then $\,^{\#}{\cal T}(\Sigma) \ge 2$.
Moreover, if $n=2$ and $\,^{\#}{\cal T}(\Sigma)<+\infty$ further
holds, then there exist at least two elliptic closed
characteristics on $\Sigma$. }

  Here the crucial point is to prove $i(x,1)\ge n$ when $(\tau,x)$
is the only geometrically distinct closed characteristic on $\Sigma$.
This conclusion is proved by using our index iteration theory and an
identity of non-degenerate closed characteristics on $\Sigma$ proved
by C. Viterbo in 1989.

  Because of Theorem 9 and other indications, we suspect that the
following holds:
\begin{equation}
 \{\,^{\#}{\cal T}(\Sigma)\,|\,\Sigma\in {\cal H}(2n)\}
  = \{k\in{\bf Z}\,|\, [\frac{n}{2}]+1\le k\le n\}\cup\{+\infty\}.
\label{22}
\end{equation}
We also suspect that closed orbits of the Reeb field on a compact
contact hypersurfaces in a symplectic manifold may have similar
properties.

  Many other problems related to iterations of periodic solution orbits
are still open, for example, the Seifert conjecture on the existence
of at least $n$ brake orbits for the given energy problem of
classical Hamiltonian systems on ${\bf R}^n$ (cf. \cite{Sei1},
\cite{ABL1} and the references there in), and the conjecture on the
existence of infinitely many geometrically distinct closed geodesics
on every compact Riemannian manifold (cf. \cite{Ban1} and the
solution for $S^2$ by J. Franks and V. Bangert). We believe that our
index iteration theory for symplectic paths and the methods we
developed to establish and apply it to nonlinear problems will have
the potential to play more roles in the study on these problems and
in other mathematical areas.

\medskip

  {\bf Acknowledgements.} The author sincerely thanks the 973 Program
of MOST, NNSF, MCME, RFDP, PMC Key Lab of MOE of China, S. S. Chern
Foundation, CEC of Tianjin, and Qiu Shi Sci. Tech. Foundation of Hong
Kong for their supports in recent years.

\label{lastpage}

\end{document}